\documentclass[12pt]{amsart}
\usepackage{amssymb,amsmath}
\newtheorem{proposition}{Proposition}[section]
\newtheorem{theorem}[proposition]{Theorem}
\newtheorem{corollary}[proposition]{Corollary}
\newtheorem{lemma}[proposition]{Lemma}

\theoremstyle{definition}

\newtheorem{remark}[proposition]{Remark}
\numberwithin{equation}{section}

\def\eps{\varepsilon}
\def\teta{\vartheta}

\def\R {\mathbb{R}}
\def\N {\mathbb{N}}
\def\H {{\mathcal H}}
\def\V {{\mathcal V}}

\def\D {{\mathcal D}}

\def\A {{\mathcal A}}

\def\F {{\mathcal F}}

\def\W {{\mathcal W}}

\def\l {\langle}
\def\r {\rangle}

\def\and{\qquad\text{and}\qquad}

\def \au {\rm}
\def \ti {\it}
\def \jou {\rm}
\def \bk {\it}
\def \no#1#2#3 {{\bf #1} (#3), #2.}
\def \eds#1#2#3 {#1, #2, #3.}

\title[Asymptotic behavior of a nonisothermal Cahn-Hilliard equation]
{Asymptotic behavior of a nonisothermal \\
viscous Cahn-Hilliard equation \\
with inertial term}

\date{\today}

\author[M.~Grasselli, H.~Petzeltov\'a, G.~Schimperna]
{Maurizio Grasselli$^\dag$, Hana Petzeltov\'a$^\ddag$, Giulio Schimperna$^\natural$}
\thanks{The first and third authors were also supported by the Italian MIUR PRIN Research
Project {\it Mo\-del\-liz\-za\-zio\-ne Matematica ed Analisi dei
Problemi a Frontiera Libera}. The second author was supported by
the Academy of Sciences of the Czech Republic, Institutional
Research Plan No. AV0Z10190503 and by Grant IAA100190606 of GA AV
CR. The third author was also supported by the HYKE Research
Training Network.}

\address{$^\dag$Dipartimento di Matematica ``F.~Brioschi''
\newline\indent
Politecnico di Milano
\newline\indent
Via Bonardi, 9
\newline\indent
I-20133 Milano, Italy}
\email{maugra@mate.polimi.it}

\address{$^\ddag$Mathematical Institute AS CR
\newline\indent
\v{Z}itn\'a, 25
\newline\indent
CZ-115 67 Praha, Czech Republic}
\email{petzelt@math.cas.cz}

\address{$^\natural$Dipartimento di Matematica ``F.~Casorati''
\newline\indent
Universit\`a degli Studi di Pavia
\newline\indent
Via Ferrata, 1
\newline\indent
I-27100 Pavia, Italy}
\email{giusch04@unipv.it}



\begin{document}

\begin{abstract}
We consider a differential model describing nonisothermal fast
phase separation processes taking place in a three-dimensional
bounded domain. This model consists of a viscous Cahn-Hilliard
equation characterized by the presence of an inertial term
$\chi_{tt}$, $\chi$ being the order parameter, which is linearly
coupled with an evolution equation for the (relative) temperature
$\teta$. The latter can be of hyperbolic type if the
Cattaneo-Maxwell heat conduction law is assumed. The state
variables and the chemical potential are subject to the
homogeneous Neumann boundary conditions. We first provide
conditions which ensure the well-posedness of the initial and
boundary value problem. Then, we prove that the corresponding
dynamical system is dissipative and possesses a global attractor.
Moreover, assuming that the nonlinear potential is real analytic,
we establish that each trajectory converges to a single steady
state by using a suitable version of the {\L}ojasiewicz-Simon
inequality. We also obtain an estimate of the decay rate to
equilibrium.
\end{abstract}

\maketitle

\section{Introduction}

\noindent
Consider a bounded domain $\Omega \subset \R^3$ with smooth boundary $\partial\Omega$ which
contains, for any time $t\ge 0$, a two-phase system subject to nonisothermal
phase separation. A well-known evolution system which describes this kind of process is
(see \cite{Cag}, cf. also \cite{BS})
\begin{equation}
\label{CPF}
\begin{cases}
&(\teta + \chi)_t - \Delta\teta = 0,\\
&\chi_t -\Delta(-\Delta\chi + \phi(\chi) -  \teta )= 0,
\end{cases}
\end{equation}
in $\Omega\times(0,\infty)$. Here $\teta$ denotes the (relative) temperature around a given critical one,
$\chi$ represents the order parameter (or phase-field) and $\phi$ is the derivative of a suitable smooth
double well potential (e.g., $\phi(r)=r^3-ar$, $a>0$). For the sake of simplicity,
all the constants have been set equal to one.

In the isothermal case,
the following singular perturbation of Cahn-Hilliard equation has been examined in several
papers (see \cite{Bo,De,GGMP1,GGMP2,ZM1,ZM2} and references therein)
\begin{equation}
\label{RCH}
\eps\chi_{tt} + \chi_t -\Delta(-\Delta\chi + \alpha\chi_t + \phi(\chi))= 0,
\end{equation}
where $\eps>0$ is a small inertial parameter and $\alpha\ge 0$ is
a viscosity coefficient. The inertial term $\eps\chi_{tt}$
accounts for fast phase separation processes (see, e.g.,
\cite{GJ}), while the motivations for introducing the viscous term
$\alpha\chi_t$ are detailed in \cite{N-C}. The above quoted works
are concerned with the analysis of the infinite-dimensional
dissipative dynamical system generated by \eqref{RCH} endowed with
suitable boundary conditions. We recall that the case $\alpha=0$
has been analyzed so far in one spatial dimension only, since in
two and three dimensions, uniqueness and smoothness of solutions
are still open issues (see however \cite{S}).

In this paper we consider equation \eqref{RCH} in the nonisothermal case, namely,
\begin{equation}
\label{CRPF}
\begin{cases}
&(\teta + \chi)_t + \nabla\cdot {\mathbf q}= 0,\\
&\sigma {\mathbf q}_t + {\mathbf q} = - \nabla\teta,\\
&\eps\chi_{tt} + \chi_t -\Delta(-\Delta\chi + \alpha\chi_t + \phi(\chi)-\teta)= 0,
\end{cases}
\end{equation}
where $\sigma\in [0,1]$. Observe that the standard Fourier law is obtained when $\sigma=0$.
Otherwise, we have the so-called Maxwell-Cattaneo heat conduction law which entails
that $\teta$ propagates at finite speed (see, e.g., \cite{HP, JP1,JP2} and their references).

System \eqref{CRPF} is subject to the initial conditions
\begin{equation}
\label{IC}
\teta(0)=\teta_0,\quad \sigma{\mathbf q}(0) = \sigma{\mathbf q}_0,\quad
\chi(0)= \chi_0,\quad
\chi_t(0)=\chi_1,
\qquad\text{ in }\Omega,
\end{equation}
and to the no-flux boundary conditions
\begin{equation}
\label{BC} {\mathbf q} \cdot {\mathbf n} = \nabla\chi\cdot{\mathbf
n}  = \nabla(\Delta\chi)\cdot{\mathbf n} =0, \qquad\text{ on
}\partial\Omega\times(0,\infty),
\end{equation}
where ${\mathbf n}$ stands for the outward normal derivative and $\cdot$
indicates the usual Euclidean scalar product.  Observe that
\eqref{CRPF} reduces to \eqref{CPF} when $\eps=\alpha=0$. Moreover, note that \eqref{BC}
are equivalent to assume the first two conditions and $\nabla u\cdot{\mathbf n}=0$, where
$u=-\Delta\chi + \alpha\chi_t + \phi(\chi)-\teta$ is the so-called chemical potential.

Here we want to demonstrate first that problem
\eqref{CRPF}-\eqref{BC} is well posed. Thus we can construct a
strongly continuous semigroup $S_\sigma(t)$ on an appropriate
phase-space. This semigroup possesses a bounded absorbing set
which is compact in the phase-space if $\sigma=0$, otherwise we
show the existence of a compact exponentially attracting set which
entails the asymptotic compactness of $S_\sigma(t)$. The latter
result is based on a recent decomposition of the solution
semigroup devised in \cite{PZ}. Therefore, for any $\sigma\geq 0$,
we deduce that $S_\sigma(t)$ possesses a (smooth) global
attractor. Taking advantage of these results, we can also deduce
that any trajectory originating from the phase-space is
precompact. Then, we can proceed to analyze the asymptotic
behavior of a single trajectory. More precisely, we show that if
$\phi$ is real analytic, then any (weak) solution
$(\teta(t),\sigma{\mathbf q}(t),\chi(t))$ converges, as $t$ goes
to $\infty$, to a single equilibrium, namely, to a triplet
$(\teta_\infty,0,\chi_\infty)$, where $\teta_\infty$ and
$\chi_\infty$ satisfy
\begin{equation}
\label{STAT}
\begin{cases}
&\teta_\infty = \vert\Omega\vert^{-1}\displaystyle\int_\Omega (\teta_0-\eps\chi_1),\\
&\displaystyle\int_\Omega \chi_\infty = \displaystyle\int_\Omega (\eps\chi_1 + \chi_0),\\
& -\Delta(-\Delta\chi_\infty + \phi(\chi_\infty))= 0,
\qquad\text{ in }\Omega,\\
&\nabla\chi_\infty\cdot{\mathbf n}  = \nabla(\Delta\chi_\infty)\cdot{\mathbf n} =0,
\qquad\text{ on }\partial\Omega.
\end{cases}
\end{equation}
This result is obtained by exploiting a well-known technique
originated from some works of S.~{\L}ojasiewicz \cite{Lo1,Lo2} and
then refined by L.~Simon \cite{Sim}. We recall that, in more than
one spatial dimension, the structure of the set of solutions to
\eqref{STAT} may contain a continuum of solutions if $\Omega$ is a
ball or an annulus (cf., e.g., \cite{Har} and references therein).
If this is the case, it is nontrivial to decide whether or not a
given trajectory converges to a single stationary state. Moreover,
this might not happen even for finite-dimensional dynamical
systems (cf. \cite{Aul}) and there are negative results for
semilinear parabolic equations with smooth nonlinearities (see
\cite{PR,PS}).

During the last years, the {\L}ojasiewicz-Simon technique has been modified and used
by many authors (cf., e.g., \cite{Chi,CF,CJ,FS,HJ,HJK,HT,Je1,Je2,WZ2}) to investigate
a number of parabolic and hyperbolic semilinear equations with variational structure.
More recently, this technique has also been used for problems with only a partial variational structure,
like the phase-field systems. More precisely, nonconserved models (with or without memory effects)
have been analyzed in \cite{AF,AFIR,FSc,GPS1,Zh}, while the case of a hyperbolic dynamics for the
order parameter has been examined in \cite{GPS2,WGZ}. There are also results for nonlocal models
(see \cite{FIP2,GPS3}). Concerning the standard Cahn-Hilliard equation, convergence to stationary states has been
examined in \cite{CFP,GG,PRZ,RH,WZ1}, while the nonconstant temperature case, namely \eqref{CPF} with \eqref{BC},
has been first analyzed in \cite{FIP1} and then in \cite{PW}
in the case of dynamic boundary conditions. The memory effects in the heat flux have been treated in
\cite{AP,AP2} for the Coleman-Gurtin law and, recently, in \cite{M} for
a generalization of the Maxwell-Cattaneo law. As we shall see, here we need a particular
{\L}ojasiewicz-Simon type inequality which is a refinement of the one proved in \cite{GG}
(see Lemma~\ref{LSINE} and its proof in Appendix).

\section{Well-posedness and uniform bounds}

\noindent
Let $H=L^2(\Omega)$ and ${\mathbf H}=(L^2(\Omega))^3$.
These spaces are endowed with the natural inner product $\l\cdot,\cdot\r$ and
the induced norm $\|\cdot\|$. For the sake of
simplicity, we will assume $\vert\Omega\vert=1$ and $\eps=1$.
Then, we set $V=H^1(\Omega)$, ${\mathbf V}=(H^1(\Omega))^3$ and $W=H^2(\Omega)$,
both endowed with their standard inner products, and we
define the subspace of $H$ of the null mean functions
$$
H_0 = \{v\in H \,:\, \l v,1\r=0 \}.
$$
We also introduce the linear nonnegative operator $A=-\Delta:\D(A)\subset H\to H_0$ with domain
$$
\D(A)= \{v\in W\,:\, \nabla v \cdot{{\mathbf n}}= 0,\quad\text{ on }\partial\Omega\},
$$
and denote by $A_0$ its restriction to $H_0$. Note that $A_0$ is a positive linear operator;
hence, for any $r\in\R$, we can define its powers $A^r$ and, consequently, set
$V_0^r=\D(A_0^{r/2})$ endowed with the inner product
$$
\langle v_1,v_2\rangle_ {V_0^r}=\langle A_0^{r/2}v_1,A_0^{r/2}v_2\rangle.
$$
Clearly, we have $V^0_0\equiv H_0$. In addition, we need to use the Hilbert spaces
$$
{\mathbf V}_0 = \{{\mathbf v}\in {\mathbf V}\,:\, {\mathbf v}\cdot{\mathbf n}=0
\quad\text{ on }\partial\Omega\},
$$
and
$$
\H_{\sigma} = H \times {\mathbf H} \times V \times V^*,
\qquad \V_{\sigma} = V \times {\mathbf V}_0 \times \D(A) \times H,
$$
endowed with the following norms, respectively,
\begin{align*}
&\Vert (z^1,{\mathbf z}^2,z^3,z^4)\Vert^2_{\H_{\sigma}}
= \Vert z^1\Vert^2 + \sigma \Vert {\mathbf z}^2\Vert^2
+ \Vert z^3\Vert^2_{V} + \Vert z^4\Vert^2_{V^*},\\
&\Vert (z^1,{\mathbf z}^2,z^3,z^4)\Vert^2_{\V_{\sigma}}
= \Vert z^1\Vert^2_V + \sigma \Vert {\mathbf z}^2\Vert^2_{{\mathbf V}}
+ \Vert z^3\Vert^2_{W} + \Vert z^4\Vert^2,
\end{align*}
if $\sigma>0$. Otherwise, we simply set
$$
\H_{0} = H \times V \times V^*,\qquad
\qquad \V_{0} = V \times \D(A) \times H.
$$

Our assumptions on the function $\phi$ and on the potential $\Phi$, defined by
$$
\Phi(y) = \int_0^y\,\phi(\xi)d\xi,\qquad\forall\,y\in\R,
$$
are the following
\begin{align}
\label{fi1}
&\Phi\in C^3(\R)\;\textrm{ such that }\; \Phi(y)\geq -c_0,\quad \forall \,y\in \R;\\
\label{fi2}
&|\phi''(y)|\leq c_1(1+|y|),\quad \forall \,y\in \R;\\
\label{fi3}
&\forall\,\epsilon>0, \;\textrm { there exists }\, c_\epsilon > 0 \;\textrm{ such that }\\
\nonumber
&\vert\phi(y)\vert \leq \epsilon \Phi(y) + c_\epsilon,\quad \forall \,y\in \R;\\
\label{fi4}
&\forall\,\varsigma\in\R, \;\textrm { there exist }\, c_2>0 \textrm{ and }\, c_3\geq 0 \;\textrm{ such that }\\
\nonumber
&(y-\varsigma)\phi(y) \geq c_2 \Phi(y) - c_3, \quad \forall \,y\in \R;\\
\label{fi5}
&\phi'(y) \geq -c_4,\quad \forall \,y\in \R;
\end{align}
for some positive constants $c_0$, $c_1$, $c_4$.
Here $c_2$ and $c_3$ continuously depend on $\varsigma$.

We now rewrite system \eqref{CRPF} together with \eqref{BC} in the following form
\begin{equation}
\label{PHabs}
\begin{cases}
&\langle (\teta + \chi)_t,v\rangle - \langle {\mathbf q},\nabla v\rangle  = 0,
\qquad\text{ in }(0,\infty),\\
&\langle \sigma{\mathbf q}_t + {\mathbf q},{\mathbf v} \rangle
= \langle \teta,\nabla\cdot {\mathbf v}\rangle,
\qquad\text{ in }(0,\infty),\\
&\langle \chi_{tt} + \chi_t,w\rangle
+ \langle A\chi + \phi(\chi) + \alpha\chi_t - \teta,Aw\rangle = 0,
\qquad\text{ in }(0,\infty),
\end{cases}
\end{equation}
for all $v\in V$, ${\mathbf v}\in {\mathbf V}_0$, and $w\in D(A)$,
endowed with initial conditions \eqref{IC}.

Let us prove
\begin{theorem}
\label{WP}
Let \eqref{fi1}-\eqref{fi5} hold. Then, for any $(\teta_0,{\mathbf q}_0,\chi_0,\chi_1)$
such that
\begin{align}
\label{I1}
&\teta_0 \in H,\\
\label{I1bis}
&\sigma{\mathbf q}_0 \in {\mathbf H},\\
\label{I2}
&\chi_0\in V,\\
\label{I3}
&\chi_1 \in V^*,
\end{align}
the Cauchy problem \eqref{PHabs}-\eqref{IC} has a (weak) solution $(\theta,\chi)$
with the following properties
\begin{align}
\label{R1}
&\teta \in C^0([0,\infty),H)\\
\label{R2}
&\sigma{\mathbf q}\in C^0([0,\infty);{\mathbf H}),\quad
{\mathbf q}\in L^2(0,\infty;{\mathbf H}),\\
\label{R3}
&\chi \in C^0([0,\infty),V),\\
\label{R4}
&\chi_t \in C^0([0,\infty),V^*) \cap L^2(0,\infty,V^*) ,\\
\label{R5}
&\alpha\chi_t \in L^2(0,\infty,H).
\end{align}
and there exists a positive constant $C$, depending on the norms of the initial data and on $\phi$,
such that, for all $t\ge 0$,
\begin{align}
\label{BOUND1}
&\Vert (\teta(t),{\mathbf q}(t),\chi(t),\chi_t(t))\Vert^2_{\H_\sigma}\\
\nonumber
&+ \int_t^\infty\,\left(\Vert \teta(\tau) - \langle\teta(\tau),1\rangle\Vert^2
+ \Vert {\mathbf q}(\tau)\Vert^2 + \Vert\chi_t(\tau)\Vert_{V^*}^2
+ \alpha\Vert\chi_t(\tau)\Vert^2\right)d\tau \leq C,
\end{align}
and
\begin{equation}
\label{MC}
\langle(\teta+\chi)(t),1\rangle = \langle \teta_0+\chi_0,1\rangle,
\qquad
\langle \chi(t),1\rangle = \langle \chi_0 + \chi_1,1\rangle - \langle \chi_1,e^{-t}\rangle.
\end{equation}
If $\alpha>0$, then the solution is unique and the following bound holds
\begin{equation}
\label{BOUND2}
\sup_{t\geq 0}\int_t^{t+1} \Vert A\chi(\tau)\Vert^2 d\tau \leq C.
\end{equation}
Moreover, for any fixed $T>0$,
if $(\teta_{0i},{\mathbf q}_{0i},\chi_{0i},\chi_{1i})\in \H_\sigma$, $i=1,2$,
then the corresponding solutions $(\teta^i,{\mathbf q}^i,\chi^i,\chi^i_t)$
satisfy
\begin{align}
\label{contdep}
&\Vert ((\teta^1-\teta^2)(t),({\mathbf q}^1-{\mathbf q}^2)(t),
(\chi^1-\chi^2)(t),(\chi^1-\chi^2)_t(t))\Vert^2_{\H_\sigma}\\
\nonumber
&\leq C(R)e^{KT}
\Vert (\teta_{01}-\teta_{02},{\mathbf q}_{01} - {\mathbf q}_{02},\chi_{01}-\chi_{02},
\chi_{11} - \chi_{12})\Vert^2_{\H_\sigma},\qquad\forall\,t\in [0,T],
\end{align}
for some positive constants $C(R)$ and $K$, both independent of $T$, where
$$
\Vert(\teta_{0i},{\mathbf q}_{0i},\chi_{0i},\chi_{1i})\Vert_{\H_\sigma}\leq R,\qquad i=1,2.
$$
\end{theorem}

\begin{proof}
We first show inequality \eqref{BOUND1} arguing formally. This argument can be made rigorous within a Faedo-Galerkin
scheme and it suffices to prove the existence of a solution for all $\alpha\geq 0$. From now on $C$ will
denote a generic positive constant which depends on $\phi$ and on
the spatial averages of the initial data, at most.
If a solution exists, then it is
easy to show the validity of \eqref{MC}, due to the boundary conditions \eqref{BC}. Moreover, we have
\begin{equation}
\label{MC1}
\langle \chi_t(t),1\rangle = \langle \chi_1,1\rangle e^{-t}.
\end{equation}
Let us set now
\begin{equation}
\label{newvar1}
\tilde \teta = \teta - \langle\teta,1\rangle, \qquad
\tilde \chi = \chi - \langle\chi,1\rangle,
\end{equation}
and rewrite problem \eqref{PHabs} in the form
\begin{equation}
\label{PHabsav}
\begin{cases}
&\langle (\tilde\teta + \tilde\chi)_t,v\rangle - \langle {\mathbf q},\nabla v\rangle  = 0,
\qquad\text{ in }(0,\infty),\\
&\langle \sigma{\mathbf q}_t + {\mathbf q},{\mathbf v} \rangle
= \langle \tilde\teta,\nabla\cdot {\mathbf v}\rangle,
\qquad\text{ in }(0,\infty),\\
&\langle \tilde\chi_{tt} + \tilde\chi_t,w\rangle
+ \langle A\tilde\chi + \phi(\chi) + \alpha\tilde\chi_t - \tilde\teta,Aw\rangle = 0,
\qquad\text{ in }(0,\infty),
\end{cases}
\end{equation}
for all $v\in V$, ${\mathbf v}\in {\mathbf V}_0$, and $w\in D(A)$.

Let us take $v=\tilde\teta$ in the first equation,
${\mathbf v}={\mathbf q}$ in the second equation, and $w=A_0^{-1}(\tilde \chi_t +
\beta \tilde \chi)$, where $\beta>0$ will be chosen small enough.
Adding together the resulting identities, we get
\begin{align}
\label{unifest1}
&\frac{d}{dt} \Big(\Vert \tilde\teta\Vert^2 + \sigma\Vert{\mathbf q}\Vert^2
+ \Vert A^{-1/2}_0 \tilde \chi_t\Vert ^2
+ \Vert \nabla \tilde \chi \Vert ^2
+ 2\beta \langle A^{-1/2}_0 \tilde\chi_t,A_0^{-1/2}\tilde\chi\rangle \\
\nonumber
&+ \beta \Vert A^{-1/2}_0 \tilde \chi \Vert ^2
+ \alpha\beta \Vert \tilde \chi \Vert ^2
+ 2\langle \Phi(\chi),1\rangle\Big)\\
\nonumber
&+ 2\Vert{\mathbf q}\Vert^2
+ 2(1-\beta)\Vert A^{-1/2}_0\tilde\chi_t\Vert^2 + 2\alpha\Vert \tilde \chi_t\Vert ^2
- 2\langle \phi(\chi),\langle \chi_t(t),1\rangle \rangle\\
\nonumber
&+ 2\beta \Vert \nabla\tilde \chi \Vert ^2
+ 2\beta\langle \phi(\chi),\tilde\chi\rangle
- 2\beta\langle\tilde\teta,\tilde\chi\rangle = 0.
\end{align}

Observe that, using \eqref{fi4} with $\varsigma=\langle \chi,1\rangle$, we deduce
\begin{equation}
\label{nonlin1}
\langle \phi(\chi),\tilde\chi\rangle \geq
C_1\langle \Phi(\chi),1\rangle - C_2,
\end{equation}
for some $C_1>0$, while, on account of \eqref{fi3}, we infer
\begin{equation}
\label{nonlin2}
-\langle \phi(\chi),\langle \chi_t,1\rangle\rangle = -\langle \phi(\chi),1\rangle \langle \chi_t,1\rangle
\geq - \frac{\beta C_1}{2}\langle \Phi(\chi),1\rangle - \frac{c_\beta}{2} e^{-t} .
\end{equation}
Hence, using \eqref{fi1}, we have
\begin{align}
\label{nonlin3}
&- 2\langle \phi(\chi),\langle \chi_t,1\rangle \rangle + 2\beta\langle \phi(\chi),\tilde\chi\rangle\\
\nonumber
&\geq \beta C_1\langle \Phi(\chi),1\rangle - 2\beta C_2 - c_\beta e^{-t}
\geq -C(\beta + c_\beta e^{-t}).
\end{align}
Then, taking \eqref{MC} and \eqref{MC1} into account, from \eqref{unifest1}
we deduce
\begin{align}
\label{unifest2}
&\frac{d}{dt} \Big(\Vert \tilde\teta\Vert^2 + \sigma\Vert{\mathbf q}\Vert^2
+ \Vert A^{-1/2}_0 \tilde \chi_t\Vert ^2
+ \Vert \nabla \tilde \chi \Vert ^2
+ 2\beta \langle A^{-1/2}_0 \tilde\chi_t,A_0^{-1/2}\tilde\chi\rangle \\
\nonumber
&+ \beta \Vert A^{-1/2}_0 \tilde \chi \Vert ^2
+ \alpha\beta \Vert \tilde \chi \Vert ^2
+ 2\langle \Phi(\chi),1\rangle\Big)\\
\nonumber
&+ 2\Vert{\mathbf q}\Vert^2
+ 2(1-\beta)\Vert A^{-1/2}_0\tilde\chi_t\Vert ^2
+ 2\alpha\Vert \tilde \chi_t\Vert ^2
+ 2\beta\Vert \nabla\tilde \chi \Vert ^2
- 2\beta\langle\tilde\teta,\tilde\chi\rangle\\
\nonumber
&\leq C(\beta + c_\beta e^{-t}).
\end{align}

Let us now test the third equation of \eqref{PHabsav} with $\tilde\chi$. We obtain
\begin{align}
\label{unifest3}
&\frac{d}{dt} \Big(2\langle \tilde\chi_{t},\tilde\chi\rangle + \Vert\tilde\chi\Vert^2
+ \alpha\Vert \nabla\tilde\chi\Vert^2\Big)\\
\nonumber
&-2\Vert\tilde\chi_t\Vert^2 + 2\Vert A\tilde\chi\Vert^2
+ 2\langle \phi'(\chi)\nabla\tilde\chi,\nabla\tilde\chi\rangle
- 2\langle \tilde\teta,A\tilde\chi\rangle = 0.
\end{align}
Moreover, in the case $\sigma>0$, using the first two equations of \eqref{PHabsav}, we have
\begin{align}
\label{unifest4}
\frac{d}{dt} \langle {\mathbf q},\nabla A^{-1}_0\tilde\teta\rangle
&= \langle {\mathbf q}_t,\nabla A^{-1}_0\tilde\teta\rangle
+ \langle {\mathbf q},\nabla A^{-1}_0\tilde\teta_t\rangle\\
\nonumber
&= -\sigma^{-1} \langle {\mathbf q},\nabla A^{-1}_0\tilde\teta\rangle
+ \sigma^{-1} \Vert \tilde\teta\Vert^2\\
\nonumber
&-\langle {\mathbf q},\nabla A^{-1}_0\tilde\chi_t\rangle
+\Vert A^{-1/2}_0\nabla\cdot{\mathbf q}\Vert^2.
\end{align}

Let us discuss first the case $\sigma>0$. Then,
multiply \eqref{unifest3} by $\gamma_1$
and \eqref{unifest4} by $-\gamma_2$, $\gamma_1>0$ and $\gamma_2>0$
to be chosen later, and sum both the obtained expressions to
\eqref{unifest2}. Note also that, by the Poincar\'e inequality and
\eqref{fi5}, for some $\kappa_1>0$ depending only on $\Omega$, we have
\begin{align}
\label{new1}
& - 2\beta\langle\tilde\teta,\tilde\chi\rangle
+ 2\gamma_1\langle \phi'(\chi)\nabla\tilde\chi,\nabla\tilde\chi\rangle
- 2\gamma_1\langle \tilde\teta,A\tilde\chi\rangle\\
\nonumber
& \ge -(\beta+2\gamma_1c_4)\Vert \nabla\tilde \chi \Vert ^2
-\gamma_1\Vert A\tilde\chi\Vert^2
- (\beta \kappa_1+\gamma_1)\Vert \tilde \teta \Vert ^2.
\end{align}
Additionally, for some $\kappa_2>0$
depending also only on $\Omega$, we get
\begin{equation}
\label{new2}
-\sigma^{-1} \langle {\mathbf q},\nabla A^{-1}_0\tilde\teta\rangle
-\langle {\mathbf q},\nabla A^{-1}_0\tilde\chi_t\rangle
\ge -\frac{\sigma^{-1}}2\Vert \tilde\teta\Vert^2
-\Vert  A_0^{-1/2}\tilde\chi_t\Vert^2
-\kappa_2\big(1+\sigma^{-1}\big)\Vert {\mathbf q}\Vert^2.
\end{equation}

Then, let us introduce the functional
\begin{align}
\label{func0}
\Psi_\sigma(\tilde\teta,{\mathbf q},\tilde\chi,\tilde\chi_t) &=\Vert \tilde\teta\Vert^2
+ \sigma\Vert{\mathbf q}\Vert^2 + \Vert A^{-1/2}_0 \tilde \chi_t\Vert ^2
+ \Vert \nabla \tilde \chi \Vert ^2\\
\nonumber
&+ 2\beta \langle A^{-1/2}_0 \tilde\chi_t,A_0^{-1/2}\tilde\chi\rangle
+ \beta \Vert A^{-1/2}_0 \tilde \chi \Vert ^2
+ \alpha\beta \Vert \tilde \chi \Vert ^2
+ 2\langle \Phi(\chi),1\rangle\\
\nonumber
&+\gamma_1\Big(2\langle \tilde\chi_{t},\tilde\chi\rangle + \Vert\tilde\chi\Vert^2
+ \alpha\Vert \nabla\tilde\chi\Vert^2\Big)\\
\nonumber
&-\gamma_2 \langle {\mathbf q},\nabla A^{-1}_0\tilde\teta\rangle,
\end{align}
and, recalling \eqref{new1} and \eqref{new2},
let us choose, in turn, $\gamma_2$ so small that
$$
\max\{2\gamma_2,\gamma_2\kappa_2(1+\sigma^{-1})\}\le1,
$$
and then $\beta$ and $\gamma_1$ so small that $\beta\le1/2$,
$\gamma_1c_4\le\beta/4$, and
$(\beta \kappa_1+\gamma_1)\le \gamma_2\sigma^{-1}/4$.
Then, $\Psi_\sigma$ fulfills the inequality
\begin{align}
\label{func1}
&\frac{d}{dt}\Psi_\sigma(\tilde\teta,{\mathbf q},\tilde\chi,\tilde\chi_t) \\
\nonumber
&+ c \Big(\sigma^{-1}\gamma_2\Vert\tilde\teta\Vert^2 + \Vert {\mathbf q}\Vert^2
+ \Vert A^{-1/2}_0\tilde\chi_t\Vert ^2
+ \alpha\Vert \tilde \chi_t\Vert ^2
+ \beta\Vert \nabla\tilde \chi \Vert ^2 \Big)  + \gamma_1\Vert A\tilde\chi\Vert^2\\
\nonumber
&\leq C(\beta + c_\beta e^{-t}).
\end{align}
Moreover, possibly choosing a smaller $\gamma_2$ (and consequently
smaller $\beta$ and $\gamma_1$), we find
\begin{align}
\label{func2}
&\Psi_\sigma(\tilde\teta(t),{\mathbf q}(t),\tilde\chi(t),\tilde\chi_t(t)) \\
\nonumber
&\geq
C_\beta\left(\Vert \teta(t) \Vert^2 + \sigma \Vert {\mathbf q}(t)\Vert^2 +
\Vert\chi(t)\Vert_{V}^2 + \Vert\chi_t(t)\Vert_{V^*}^2\right)  -C, \qquad\forall\,t\geq 0.
\end{align}
On the other hand, on account of \eqref{fi1}, \eqref{fi2} and \eqref{I1}-\eqref{I3},
and recalling notation \eqref{newvar1}, we find $R_0>0$ such that
$$
\Psi_\sigma(\tilde\teta_0,{\mathbf q}_0,\tilde\chi_0,\tilde\chi_1)\leq R_0.
$$
Using then \cite[Lemma~2.1]{GGMP1}, we deduce that there exists $t_0=t_0(R_0)>0$ such that,
for all $t\geq t_0$,
$$
\Psi_\sigma(\tilde\teta(t),{\mathbf q}(t),\tilde\chi(t),\tilde\chi_t(t)) \leq R,
$$
where $R$ is independent of $R_0$. Thus, recalling \eqref{func2}, we deduce that
\begin{equation}
\label{BOUND3}
\Vert (\teta(t),{\mathbf q}(t),\chi(t),\chi_t(t))\Vert^2_{\H_\sigma} \leq C(R_0),
\end{equation}
for all $t\in [0,\infty)$.
On account of \eqref{MC1} and \eqref{BOUND3}, taking $\beta=0$ in \eqref{unifest1},
integrating from $t$ to $T$ and letting $T$ go to $\infty$ we also get the integral control
\begin{equation}
\label{BOUND4}
\int_t^\infty\,\left(
\Vert {\mathbf q}(\tau)\Vert^2 + \Vert\chi_t(\tau)\Vert_{V^*}^2
+ \alpha\Vert\chi_t(\tau)\Vert^2\right)d\tau \leq C(R_0).
\end{equation}
Then, using \eqref{BOUND3} and \eqref{BOUND4}, from \eqref{unifest4} we deduce
$$
\int_t^\infty\, \Vert \tilde\teta(\tau)\Vert^2 d\tau \leq C(R_0),
$$
so that \eqref{BOUND1} is proved. In addition, integrating \eqref{func1} from $t$ to $t+1$,
we then find \eqref{BOUND2}.

The case $\sigma=0$ is simpler. We can take the functional
\begin{align*}
\Psi_0(\tilde\teta,\tilde\chi,\tilde\chi_t) &=\Vert \tilde\teta\Vert^2
+ \Vert A^{-1/2}_0 \tilde \chi_t\Vert ^2
+ \Vert \nabla \tilde \chi \Vert ^2\\
\nonumber
&+ 2\beta \langle A^{-1/2}_0 \tilde\chi_t,A_0^{-1/2}\tilde\chi\rangle
+ \beta \Vert A^{-1/2}_0 \tilde \chi \Vert ^2
+ \alpha\beta \Vert \tilde \chi \Vert ^2
+ 2\langle \Phi(\chi),1\rangle\\
\nonumber
&+\gamma_1\Big(2\langle \tilde\chi_{t},\tilde\chi\rangle + \Vert\tilde\chi\Vert^2
+ \alpha\Vert \nabla\tilde\chi\Vert^2\Big),
\end{align*}
and observe that
\begin{align*}
&\frac{d}{dt}\Psi_0(\tilde\teta,\tilde\chi,\tilde\chi_t)
+ c \Big(\Vert \nabla\tilde\teta\Vert^2
+ \Vert A^{-1/2}_0\tilde\chi_t\Vert^2
+ \alpha\Vert \tilde \chi_t\Vert ^2
+ \beta\Vert \nabla\tilde \chi \Vert ^2 \Big)  + \gamma_1\Vert A\tilde\chi\Vert^2\\
\nonumber
&\leq C(\beta + c_\beta e^{-t}).
\end{align*}
Then we can argue as above.

Estimate \eqref{contdep} is standard, provided that $\alpha>0$.
Indeed, it suffices to write down problem \eqref{PHabs} for the
difference of two solutions $(\teta_i,{\mathbf q}_i,\chi_i)$,
$i=1,2$, and then multiply the first equation by
$\teta_1-\teta_2$, the second one by ${\mathbf q}_1 - {\mathbf
q}_2$, and the third one by $A_0^{-1}(\tilde \chi_1-\tilde
\chi_2)_t$. Using the Gronwall lemma and taking \eqref{fi2} into
account, one easily gets the wanted estimate (see, e.g., \cite{Bo}
or \cite{GGMP1} for the isothermal case).
\end{proof}From Theorem~\ref{WP} and its proof we deduce that, letting
$$
X^\delta_\sigma= \{(z^1,{\mathbf z}^2, z^3,z^4)\in \H_\sigma\; :\;
\vert\langle z^1,1\rangle\vert + \vert\langle z^3,1\rangle\vert +
\vert\langle z^4,1\rangle\vert \leq \delta\}
$$
for some $\delta\geq 0$, endowed with the metric induced by the norm of $\H_\sigma$,
and setting
$$
(\teta(t),{\mathbf q}(t),\chi(t),\chi_t(t))=:S_\sigma(t)(\teta_0,{\mathbf q}_0,\chi_0,\chi_1), \qquad\forall\,t\geq 0,
$$
we have that $S_\sigma(t)$ is a strongly continuous semigroup on $X^\delta_\sigma$
with a bounded absorbing set. Similarly, we can define a strongly continuous dissipative
semigroup $S_0(t)$ on $X^\delta_0$. Summing up, we have

\begin{corollary}
\label{bddabs}
Let \eqref{fi1}-\eqref{fi5} hold.
For any given $\sigma\in [0,1]$, the semigroup $S_\sigma(t)$ acting on $X^\delta_\sigma$
has a bounded absorbing set.
\end{corollary}


\section{Precompactness of trajectories and global attractor}

\noindent
Here we prove

\begin{theorem}
\label{precomp}
Let \eqref{fi1}-\eqref{fi5} hold and suppose $\alpha>0$. If $\sigma\in (0,1]$ and
$(\teta_0,{\mathbf q}_0,\chi_0,\chi_1)$ satisfies
\eqref{I1}-\eqref{I3}, then, indicating by $(\teta,{\mathbf q},\chi)$
the corresponding solution to \eqref{PHabs}-\eqref{IC} given by Theorem~\ref{WP},
the orbit $\bigcup_{t\geq 0}\, (\teta(t),{\mathbf q}(t),\chi(t),\chi_t(t))$
is precompact in $\H_\sigma$. Moreover, there holds
\begin{align}
\label{C1}
&\Vert\teta(t)-\langle \teta_0 -\chi_1,1\rangle\Vert \to 0,\\
\label{C2}
&\Vert{\mathbf q}(t)\Vert \to 0,\\
\label{C3}
&\Vert \chi_t(t)\Vert_{V^*} \to 0,
\end{align}
as $t$ goes to $\infty$, and the $\omega$-limit set $\omega(\teta_0,{\mathbf q}_0,\chi_0,\chi_1)$ consists only of
equilibrium points of the form $(\teta_\infty,0,\chi_\infty,0)$ where $(\teta_\infty,\chi_\infty)$ satisfies \eqref{STAT}.
Similar results hold when $\sigma=0$.
\end{theorem}

\begin{proof}
On account of \cite{PZ}, observe first that, thanks to \eqref{fi2}, \eqref{fi5},
and \eqref{BOUND1}, we can choose $\ell\geq c_4$ large enough, and depending on the norms of the initial data,
such that
\begin{equation}
\label{ineqfi}
\frac{1}{2}\Vert \nabla z \Vert^2 + (\ell-2c_4)\Vert z\Vert^2
- \langle \phi'(\chi(t))z,z\rangle \geq 0,
\end{equation}
for all $z\in V$ and every $t\geq 0$. Consequently, we set
$$
\psi(r) = \phi(r) + \ell r, \qquad\forall\,r\in\R.
$$
Then, we split the solution to \eqref{PHabs} in this way
$$
(\teta,{\mathbf q},\chi) = (\teta^d,{\mathbf q}^d,\chi^d) + (\teta^c,{\mathbf q}^c,\chi^c),
$$
where
\begin{equation}
\label{PHabs1}
\begin{cases}
&\langle (\teta^d+ \chi^d)_t,v\rangle - \langle {\mathbf q}^d,\nabla v\rangle  = 0,
\qquad\text{ in }(0,\infty),\\
&\langle \sigma{\mathbf q^d}_t + {\mathbf q^d},{\mathbf v} \rangle
= \langle \teta_d,\nabla\cdot {\mathbf v}\rangle,
\qquad\text{ in }(0,\infty),\\
&\langle \chi^d_{tt} + \chi^d_t,w\rangle
+ \langle A\chi^d + \psi(\chi) - \psi(\chi^c) + \alpha\chi^d_t - \teta^d,Aw\rangle = 0,
\qquad\text{ in }(0,\infty),\\
&\teta^d(0)=\tilde\teta_0,\quad
\sigma{\mathbf q}^d(0) = \sigma{\mathbf q}_0,\quad
\chi^d(0)= \tilde \chi_0,\quad
\chi^d_t(0)=\tilde \chi_1,
\qquad\text{ in }\Omega,
\end{cases}
\end{equation}
and
\begin{equation}
\label{PHabs2}
\begin{cases}
&\langle (\teta^c+ \chi^c)_t,v\rangle - \langle {\mathbf q}^c,\nabla v\rangle  = 0,
\qquad\text{ in }(0,\infty),\\
&\langle \sigma{\mathbf q^c}_t + {\mathbf q^c},{\mathbf v} \rangle
= \langle \teta_c,\nabla\cdot {\mathbf v}\rangle,
\qquad\text{ in }(0,\infty),\\
&\langle \chi^c_{tt} + \chi^c_t,w\rangle
+ \langle A\chi^c + \psi(\chi^c) + \alpha\chi^c_t - \teta^c,Aw\rangle
= \langle \ell\chi,Aw\rangle,
\qquad\text{ in }(0,\infty),\\
&\teta^c(0)=\langle\teta_0,1\rangle,\quad
\sigma{\mathbf q}^c(0) = {\mathbf 0},\quad
\chi^c(0)= \langle \chi_0,1\rangle,\quad
\chi^c_t(0)=\langle \chi_1,1\rangle,
\qquad\text{ in }\Omega,
\end{cases}
\end{equation}
for all $v\in V$, ${\mathbf v}\in {\mathbf V}_0$, and $w\in D(A)$.

We shall prove that
$(\teta^d(t),{\mathbf q}^d(t),\chi^d(t),\chi^d_t(t))$
exponentially decays at $0$ in $\H_\sigma$ as $t$ goes to $\infty$, while
$(\teta^c,{\mathbf q}^c,\chi^c,\chi^c_t)$ is bounded in a space which is
compactly embedded in $\H_\sigma$, uniformly in time.

Let us prove first that, for any $t\geq s \geq 0$ and every $\varpi>0$, there
holds
\begin{equation}
\label{int1}
\alpha\int_s^t \Vert \chi^c_t(\tau)\Vert^2 d\tau \leq \varpi(t-s) + \frac{C}{\varpi}.
\end{equation}
This estimate combined with \eqref{BOUND1} will allow us to use a suitable version of
the Gronwall Lemma (see \cite[Lemma~5]{PZ}).

Let us take $v=\tilde\teta^c$ in the first equation of \eqref{PHabs2},
${\mathbf v}={\mathbf q}^c$ in the second equation, and $w=A_0^{-1}\tilde\chi^c_t$
in the third one. Then we obtain
\begin{align}
\label{unifest5}
&\frac{d}{dt} \Big(\Vert \tilde\teta^c\Vert^2 + \sigma\Vert{\mathbf q}^c\Vert^2
+ \Vert A^{-1/2}_0 \tilde \chi^c_t\Vert ^2
+ \Vert \nabla \tilde \chi^c \Vert ^2
+ 2\langle \Psi(\chi^c),1\rangle
-2\ell\langle \chi,\tilde\chi^c\rangle \Big)\\
\nonumber
&+ 2\Vert{\mathbf q}^c\Vert^2
+ 2 \Vert A^{-1/2}_0\tilde\chi^c_t\Vert + 2\alpha\Vert \tilde \chi^c_t\Vert ^2\\
\nonumber
&=2\langle \psi(\chi^c),\langle \chi_t,1\rangle \rangle
- 2\ell\langle \chi_t,\tilde\chi^c\rangle.
\end{align}
Here $\Psi$ is a primitive of $\psi$. Observe first that it is not difficult to realize that
an estimate similar to  \eqref{BOUND1} holds for $(\teta^c,{\mathbf q}^c,\chi^c,\chi^c_t)$
as well. Therefore, on account of \eqref{fi2} and \eqref{MC1}, we have, for any $\varpi>0$ and any $t\geq 0$,
$$
2\langle \psi(\chi^c(t)),\langle \chi_t(t),1\rangle \rangle
- 2\ell\langle \chi_t(t),\tilde\chi^c(t)\rangle
\leq  2\varpi + \frac{C}{\varpi} (\Vert \chi_t(t)\Vert^2 + \varpi e^{-t})
$$
Therefore, \eqref{int1} follows from integrating \eqref{unifest5}
with respect to time from $s$ to $t$, using the above inequality and \eqref{BOUND1},
recalling that $\alpha>0$, and observing that (cf.~\eqref{MC1})
$$
\langle \chi^c_t(t),1\rangle = \langle \chi_1,1\rangle e^{-t}.
$$

In order to prove the exponential decay of $(\teta^d,{\mathbf q}^d,\chi^d,\chi^d_t)$,
we first note that (cf.~\eqref{MC})
\begin{equation}
\label{MC2}
\langle \teta^d(t),1\rangle = \langle \chi^d(t),1\rangle = 0, \qquad\forall\,t\geq 0,
\end{equation}
so that $\teta^d=\tilde\teta^d$ and $\chi^d=\tilde\chi^d$.

We now argue as to get \eqref{func0}, namely, we take $v=\teta^d$ in the first equation,
${\mathbf v}={\mathbf q}^d$ in the second equation, and
$w=A_0^{-1}(\chi^d_t + \beta\chi^d)$, with some $\beta>0$ which will be chosen
later on. Then, we add the functional
$-\gamma\langle {\mathbf q}^d,\nabla A^{-1}\teta^d\rangle$ with $\gamma>0$.
Thus, recalling \eqref{unifest4}, we obtain
\begin{align}
\label{decay1}
&\frac{d}{dt} \Big(\Vert \teta^d\Vert^2 + \sigma\Vert{\mathbf q}^d\Vert^2
+ \Vert A^{-1/2}_0 \chi^d_t\Vert ^2
+ \Vert \nabla \chi^d \Vert ^2
+ 2\beta \langle A^{-1/2}_0 \chi^d_t,A_0^{-1/2}\chi^d\rangle \\
\nonumber
&+ \beta \Vert A^{-1/2}_0 \chi^d \Vert ^2
+ \alpha\beta \Vert \chi^d \Vert ^2
-\gamma\langle {\mathbf q}^d,\nabla A^{-1}\teta^d\rangle\\
\nonumber
&+ 2\langle \psi(\chi)-\psi(\chi^c),\chi^d\rangle
- \langle \psi'(\chi)\chi^d,\chi^d\rangle\Big)\\
\nonumber
&+ 2\Vert{\mathbf q}^d\Vert^2
+ 2(1-\beta)\Vert A^{-1/2}_0\chi^d_t\Vert ^2 + 2\alpha\Vert\chi^d_t\Vert ^2\\
\nonumber
&-\frac{\gamma}{\sigma} \langle {\mathbf q}^d,\nabla A^{-1}_0 \teta^d\rangle
+ \frac{\gamma}{\sigma} \Vert \teta^d\Vert^2
-\gamma\langle {\mathbf q}^d,\nabla A^{-1}_0 \chi^d_t\rangle
+\gamma\Vert A^{-1/2}_0\nabla\cdot{\mathbf q}^d\Vert^2\\
\nonumber
&+ 2\beta \Vert \nabla \chi^d \Vert ^2
+ 2\beta\langle \psi(\chi)-\psi(\chi^c),\chi^d\rangle
- 2\beta\langle\teta^d,A_0^{-1}\chi^d\rangle \\
\nonumber
&= 2\langle (\psi'(\chi) - \psi'(\chi^c))\chi^c_t,\chi^d\rangle
-\langle \psi^{''}(\chi)\chi_t,(\chi^d)^2\rangle.
\end{align}
Observe that, owing to \eqref{fi1}, \eqref{BOUND1}, and \eqref{MC1}, we have
\begin{align}
\label{decay2}
&2\langle (\psi'(\chi) - \psi'(\chi^c))\chi^c_t,\chi^d\rangle
-\langle \psi^{''}(\chi)\chi_t,(\chi^d)^2\rangle\\
\nonumber
&\leq C\left(\Vert \chi_t\Vert + \Vert \chi^c_t\Vert\right)\Vert \nabla\chi^d\Vert^2\\
\nonumber
&\leq \beta\Vert \nabla\chi^d\Vert^2 +
C_\beta\left(\Vert \chi_t\Vert^2 + \Vert \chi^c_t\Vert^2\right)
\Vert \nabla\chi^d\Vert^2.
\end{align}
On the other hand, setting
\begin{align*}
\Lambda_d&=\Vert \teta^d\Vert^2 + \sigma\Vert{\mathbf q}^d\Vert^2
+ \Vert A^{-1/2}_0 \chi^d_t\Vert ^2
+ \Vert \nabla \chi^d \Vert ^2
+ 2\beta \langle A^{-1/2}_0 \chi^d_t,A_0^{-1/2}\chi^d\rangle \\
\nonumber
&+ \beta \Vert A^{-1/2}_0 \chi^d \Vert ^2
+ \alpha\beta \Vert \chi^d \Vert ^2 -\gamma\langle {\mathbf q}^d,\nabla A^{-1}\teta^d\rangle\\
&+ 2\langle \psi(\chi)-\psi(\chi^c),\chi^d\rangle
- \langle \psi'(\chi)\chi^d,\chi^d\rangle,
\end{align*}
and observing that (cf.~\eqref{fi5} and \eqref{ineqfi})
$$
2\langle \psi(\chi)-\psi(\chi^c),\chi^d\rangle
- \langle \psi'(\chi)\chi^d,\chi^d\rangle
\geq (\ell - 2c_4)\Vert \chi^d\Vert^2 -
\langle \phi'(\chi)\chi^d,\chi^d\rangle
\geq -\frac{1}{2}\Vert \nabla \chi^d\Vert^2,
$$
we have that, for $\beta$ and $\gamma$ small enough,
\begin{equation}
\label{noreq}
\frac{1}{4}\Vert (\teta^d,{\mathbf q}^d,\chi^d,\chi^d_t)\Vert^2_{\H_\sigma}
\leq \Lambda_d \leq C
\Vert (\teta^d,{\mathbf q}^d,\chi^d,\chi^d_t)\Vert^2_{\H_\sigma}.
\end{equation}
Moreover, possibly choosing $\beta$ and $\gamma$ smaller than before,
and using \eqref{decay2}, from \eqref{decay1} we infer
$$
\frac{d}{dt}\Lambda_d + c_{\beta,\gamma} \Lambda_d \leq
C_{\beta,\gamma}\left(\Vert \chi_t\Vert^2 + \Vert \chi^c_t\Vert^2\right)\Lambda_d.
$$
Thus, on account of \eqref{BOUND1} and \eqref{int1}, we can apply \cite[Lemma~5]{PZ}
and deduce the exponential decay of $\Lambda_d$, so that (cf.~\eqref{MC1} and \eqref{noreq})
\begin{equation}
\label{decay3}
\Vert (\teta^d(t),{\mathbf q}^d(t),\chi^d(t),\chi^d_t(t))\Vert_{\H_\sigma}
\leq C(R)e^{-ct},
\end{equation}
provided that $\Vert(\teta_0,{\mathbf q}_0,\chi_0,\chi_1)\Vert_{\H_\sigma}\leq R$.

Moreover, taking $w=\chi^d$ in the third equation of \eqref{PHabs1},
we obtain (cf.~\eqref{unifest3})
\begin{align*}
&\frac{d}{dt} \Big(2\langle\chi^d_{t},\chi^d\rangle
+ \Vert\chi^d\Vert^2
+ \alpha\Vert \nabla \chi^d\Vert^2\Big)\\
\nonumber
&-2\Vert\chi^d_t\Vert^2 + 2\Vert A\chi^d\Vert^2
+ 2\langle \psi(\chi)-\psi(\chi^c),A\chi^d\rangle
- 2\langle \teta^d,A\chi^d\rangle = 0,
\end{align*}
which yields, using the Young inequality, \eqref{unifest1}, and \eqref{decay3},
\begin{align*}
&\frac{d}{dt} \Big(2\langle\chi^d_{t},\chi^d\rangle
+ \Vert\chi^d\Vert^2
+ \alpha\Vert \nabla \chi^d\Vert^2\Big) + \Vert A\chi^d\Vert^2
\leq C(1 + \Vert\chi^d_t\Vert^2).
\end{align*}
On account of \eqref{BOUND4} and \eqref{int1}, we have
$$
\sup_{t\geq 0} \int_t^{t+1} \Vert \chi^d_t(\tau)\Vert^2 d\tau
\leq C,
$$
so that the additional bound holds
\begin{align}
\label{int2}
\sup_{t\geq 0} \int_t^{t+1} \Vert A\chi^d(\tau)\Vert^2 d\tau
\leq C.
\end{align}

We now consider \eqref{PHabs2}. Taking $v=A\teta^c$ in the first equation,
${\mathbf v}= -\nabla\nabla\cdot{\mathbf q}^c$ in the second one,
and adding together the resulting identities, we obtain
\begin{equation}
\label{smooth1}
\frac{d}{dt}\left(\Vert\nabla\teta^c\Vert^2
+ \sigma \Vert \nabla\cdot{\mathbf q}^c\Vert^2\right)
+2\langle \chi^c_t,A\teta^c\rangle  + 2\Vert \nabla\cdot{\mathbf q}^c\Vert^2=0.
\end{equation}
We also have (cf.~\eqref{unifest4})
\begin{align}
\label{smooth2}
\frac{d}{dt} \langle {\mathbf q}^c,\nabla\teta^c\rangle
&= \langle {\mathbf q}^c_t,\nabla \teta^c\rangle
+ \langle {\mathbf q}^c,\nabla\teta^c_t\rangle\\
\nonumber
&= -\sigma^{-1} \langle {\mathbf q}^c,\nabla \teta^c\rangle
+ \sigma^{-1} \Vert \nabla\teta^c\Vert^2\\
\nonumber
&-\langle {\mathbf q}^c,\nabla\chi^c_t\rangle
+\Vert \nabla\cdot{\mathbf q}^c\Vert^2.
\end{align}
Let us now take $w=\chi^c_t + \beta\chi^c$ in the third equation. We find
\begin{align}
\label{smooth3}
&\frac{d}{dt} \Big(\Vert \chi^c_t\Vert ^2
+ \Vert A \chi^c \Vert ^2
+ 2\beta \langle\chi^c_t,\chi^c\rangle \\
\nonumber
&+ \beta \Vert \chi^c \Vert ^2
+ \alpha\beta \Vert \nabla \chi^c \Vert ^2
+ 2\langle \psi(\chi^c),A\chi^c\rangle\Big)\\
\nonumber
&+ 2(1-\beta)\Vert \chi^c_t\Vert^2 + 2\alpha\Vert \nabla\chi^c_t\Vert ^2
+ 2\beta \Vert A \chi^c \Vert ^2 \\
\nonumber
&+ 2\beta\langle \psi(\chi^c),A\chi^c\rangle
-2\langle\teta^c,A\chi^c_t\rangle - 2\beta\langle\teta^c,A\chi^c\rangle\\
\nonumber
&-2\langle \psi'(\chi^c)\chi^c_t,A\chi^c\rangle = \langle \ell\nabla\chi,\nabla (\chi^c_t + \beta\chi^c)\rangle.
\end{align}
Observe that, on account of \eqref{BOUND1},
\begin{align}
\label{crucial}
\langle \psi'(\chi^c)\chi^c_t,A\chi^c\rangle
&\leq C \left(1 + \Vert\chi^c\Vert^2_{L^6(\Omega)}\right)\Vert\chi^c_t\Vert_{L^6(\Omega)}
\Vert A\chi^c\Vert\\
\nonumber
&\leq C \Vert\chi^c_t\Vert_{V}
\Vert A\chi^c\Vert.
\end{align}
Therefore, setting
\begin{align*}
\Lambda_c &= \Vert\nabla\teta^c\Vert^2
+ \sigma \Vert \nabla\cdot{\mathbf q}^c\Vert^2
+\Vert \chi^c_t\Vert ^2
+ \Vert A \chi^c \Vert ^2
+ 2\beta \langle \chi^c_t,\chi^c\rangle \\
\nonumber
&+ \beta \Vert \chi^c \Vert ^2
+ \alpha\beta \Vert \nabla \chi^c \Vert ^2
+ 2\langle \psi(\chi^c),A\chi^c\rangle
-\gamma\langle {\mathbf q}^c,\nabla\teta^c\rangle,
\end{align*}
for some $\gamma>0$, using the Young inequality,
we can choose $\beta$ and $\gamma$ small enough so that
$$
\frac{d}{dt}\Lambda_c + c_{\beta,\gamma} \Lambda_c \leq
C_{\beta,\gamma}\left(1+ \Vert A \chi^c\Vert^2\right).
$$
Then, on account of \eqref{BOUND2} and \eqref{int2},
we obtain the uniform boundedness
of $\Lambda_c$ which implies
\begin{equation}
\label{smooth4}
\Vert\nabla\teta^c(t)\Vert^2
+ \sigma \Vert \nabla\cdot{\mathbf q}^c(t)\Vert^2
+\Vert \chi^c_t(t)\Vert ^2
+ \Vert A \chi^c(t) \Vert ^2 \leq C,\qquad\forall\,t\geq 0.
\end{equation}
The second equation of \eqref{PHabs2} can now be written in the strong form,
namely,
$$
\sigma{\mathbf q^c}_t + {\mathbf q^c}=-\nabla\teta_c,
\qquad\textrm{ a.e. in }\,\Omega\times(0,\infty),
$$
so that
$$
\sigma(\nabla\times{\mathbf q^c})_t + \nabla\times{\mathbf q^c}= {\mathbf 0},
\qquad\textrm{ a.e. in }\, \Omega\times(0,\infty),
$$
and, since $(\nabla\times {\mathbf q}^c)(0)={\mathbf 0}$, we have
$(\nabla\times{\mathbf q^c})(t)={\mathbf 0}$ for any $t\geq 0$. Consequently,
thanks to \eqref{smooth4}, $\Vert {\mathbf q}^c(t)\Vert_{{\mathbf V}}$
is uniformly bounded as well.

Summing up, we have shown that a given trajectory originating from $\H_\sigma$
is a sum of an exponentially decaying part and a term which belongs to a
closed bounded subset of $\V_\sigma$.
Therefore the trajectory is precompact in $\H_\sigma$ and, due to the integral controls
of \eqref{BOUND1} and to \eqref{MC}, we infer \eqref{C1}-\eqref{C3}. Finally,
it is not difficult to prove that
$$
\omega(\teta_0,{\mathbf q}_0,\chi_0,\chi_1)\subseteq \{(\teta_\infty,0,\chi_\infty,0)\,:\,
(\teta_\infty,\chi_\infty) \;\textrm{ satisfies }\; \eqref{STAT}\}.
$$

The case $\sigma=0$ is easier. In fact, arguing as in the isothermal case (see \cite{Bo}),
we can prove the bound
$$
\Vert (\teta(t),\chi(t),\chi_t(t))\Vert^2_{\V_0}
\leq C, \qquad \forall\,t\geq t_1=t_1(R)>0,
$$
provided that $\Vert(\teta_0,\chi_0,\chi_1)\Vert_{\H_0}\leq R$. Hence the trajectory
is precompact in $\H_0$ and we can conclude as above.
\end{proof}From the proof of Theorem~\ref{precomp}, we
deduce that the semigroup $S_\sigma(t)$ has a bounded attracting set in
$\V_\sigma$, for any $\sigma\in (0,1]$, while $S_0(t)$ has a compact
absorbing set. Therefore we have (see, e.g., \cite{Ha,Te})

\begin{corollary}
\label{globattr}
For each $\sigma\in [0,1]$, the semigroup $S_\sigma(t)$ has a connected global attractor
$\A_\sigma$ which is bounded in $\V_\sigma$.
\end{corollary}

\begin{remark}
The above result is a first, but essential, step toward the construction of a family
of exponential attractors which is stable (robust) with respect to $\sigma$ and, possibly,
to $\eps$ (see \cite{GGMP1} for the isothermal case). This will be the subject of a future
investigation.
\end{remark}


\section{Convergence to stationary states}

\noindent
Let us set
$$
E(v) = \frac{1}{2}\Vert\nabla v \Vert^2 + \langle \tilde\Phi(v),1\rangle
$$
for any $v\in V^1_0$, where
$$
\widetilde\Phi(y) = \int_0^y\,\tilde\phi(\xi)d\xi,\qquad\forall\,y\in\R,
$$
and
\begin{equation}
\label{newphi}
\tilde\phi(y) = \phi(y+\langle \chi_0 + \chi_1,1\rangle),\qquad \forall\,y\in \R.
\end{equation}

The version of the {\L}ojasiewicz-Simon inequality we need is the following
(see Appendix)

\begin{lemma}
\label{LSINE}
Suppose that $\phi$ is real analytic and assume \eqref{fi2} and \eqref{fi5}.
Let $v_\infty\in V^2_0$ be such that
\begin{equation}
\label{S}
A(A_0v_\infty + \tilde\phi(v_\infty))= 0.
\end{equation}
Then there exist $\rho\in (0,\frac12)$, $\eta>0$, and
a positive constant $L$ such that
\begin{equation}
\label{LS}
\vert E(v) - E(v_\infty)\vert ^{1-\rho} \leq L
\Vert A_0v + \tilde\phi(v) - \langle \tilde\phi(v),1\rangle \Vert_{V_0^{-1}},
\end{equation}
for all $v\in V_0^1$ such that $\Vert v - v_\infty\Vert_{V_0^1}\leq \eta$.
\end{lemma}

Then we prove
\begin{theorem}
\label{CONVRATE}
Let the assumptions of Lemma~\ref{LSINE} hold and let $\alpha>0$ and $\sigma>0$ be fixed.
If $(\teta_0,{\bf q}_0,\chi_0,\chi_1)$ satisfies
\eqref{I1}-\eqref{I3}, then the trajectory $(\teta(t),{\bf q}(t),\chi(t),\chi_t(t))$ originated from
$(\teta_0,{\bf q}_0,\chi_0,\chi_1)$ is such that
\begin{equation}
\label{singleton}
\omega(\teta_0,{\bf q}_0,\chi_0,\chi_1) = \{(\teta_\infty,{\bf 0},\chi_\infty,0)\},
\end{equation}
where $(\teta_\infty,\chi_\infty)$ satisfies
\begin{equation}
\label{STAT1}
\begin{cases}
&\teta_\infty = \vert\Omega\vert^{-1}\displaystyle\int_\Omega (\teta_0-\chi_1),\\
&\displaystyle\int_\Omega \chi_\infty = \displaystyle\int_\Omega (\chi_1 + \chi_0),\\
& A(A\chi_\infty + \phi(\chi_\infty))= 0.
\end{cases}
\end{equation}
Moreover,
\begin{equation}
\label
{convchi}
\lim_{t\to \infty} \Vert\chi(t) - \chi_\infty\Vert_V=0,
\end{equation}
and there exists $t^*>0$ and a positive constant $C$ such that
\begin{equation}
\label{rate1}
\Vert\teta(t) - \teta_\infty\Vert_{V^*} +
\|\chi(t) - \chi_\infty\|_{V^*}  \le C t^{-\frac{\rho}{1-2\rho}},
\qquad \forall\,t\ge t^*.
\end{equation}
If $\sigma=0$ a similar result hold.
\end{theorem}

\begin{proof}

Let us set $\sigma=1$ for simplicity. On account of Theorem~\ref{precomp}, we consider
$$
(\teta_\infty,0,\chi_\infty,0)\in \omega(\teta_0,{\mathbf q}_0,\chi_0,\chi_1),
$$
and we observe first that \eqref{C1}-\eqref{C3} hold and $(\teta_\infty,\chi_\infty)$
fulfills \eqref{STAT1}.

On account of \eqref{MC}, we can rewrite \eqref{PHabsav} in the form
\begin{equation}
\label{PHabsav2}
\begin{cases}
&\langle (\tilde\teta + \tilde\chi)_t,v\rangle - \langle {\mathbf q},\nabla v\rangle  = 0,
\qquad\text{ in }(0,\infty),\\
&\langle \sigma{\mathbf q}_t + {\mathbf q},{\mathbf v} \rangle
= \langle \tilde\teta,\nabla\cdot {\mathbf v}\rangle,
\qquad\text{ in }(0,\infty),\\
&\langle \tilde\chi_{tt} + \tilde\chi_t,w\rangle
+ \langle A_0\tilde\chi + \tilde\phi(\tilde\chi) + \alpha\tilde\chi_t
- \tilde\teta,Aw\rangle = \langle h(\tilde\chi),Aw\rangle,
\qquad\text{ in }(0,\infty),
\end{cases}
\end{equation}
for all $v\in V$, ${\mathbf v}\in {\mathbf V}_0$, and $w\in \D(A)$. Here we have set (cf. \eqref{newphi})
\begin{equation}
\label{newvar2}
h(\tilde\chi)=\tilde\phi(\tilde\chi) - \tilde\phi(\tilde\chi - \langle \chi_1,e^{-t}\rangle),
\qquad \forall\,r\in \R, \;\forall\,t\geq 0.
\end{equation}

Arguing as in the proof of Theorem~\ref{WP} (cf.~\eqref{unifest1}), we find
\begin{equation}
\label{Lyap1}
\frac{d}{dt}{\mathcal L}(\tilde\teta(t),{\mathbf
  q}(t),\tilde\chi(t),\tilde\chi_t(t))
 = -\Vert{\mathbf q}(t)\Vert^2 - \Vert\tilde\chi_t(t)\Vert^2_{V*} -\alpha\Vert\tilde\chi_t(t)\Vert^2
+ \langle h(\tilde\chi(t)),\tilde\chi_t(t)\rangle,
\end{equation}
where
\begin{align}
\label{Lyap2}
&{\mathcal L}(\tilde\teta(t),{\mathbf q}(t),\tilde\chi(t),\tilde\chi_t(t))\\
\nonumber
&= \frac{1}{2}\left(\Vert\tilde\teta(t)\Vert^2 + \Vert {\mathbf q}(t)\Vert^2
+ \Vert\nabla\tilde\chi(t)\Vert^2
+ 2\langle \tilde\Phi(\tilde\chi(t)),1\rangle + \Vert\tilde\chi_t(t)\Vert^2_{V*}\right).
\end{align}
Note that, due to \eqref{fi1}, \eqref{fi2}, \eqref{BOUND1} and \eqref{newvar2}, there holds
$$
\langle h(\tilde\chi(t)),\tilde\chi_t(t)\rangle \leq C_\alpha e^{-2t}
+ \frac{\alpha}{2}\Vert\tilde\chi_t(t)\Vert^2,
$$
using also the Young inequality. Therefore, from \eqref{Lyap1} we deduce
\begin{equation}
\label{Lyap3}
\frac{d}{dt}{\mathcal L}(\tilde\teta(t),{\mathbf q}(t),\tilde\chi(t),\tilde\chi_t(t)) \leq
-\Vert{\mathbf q}(t)\Vert^2 - \Vert\tilde\chi_t(t)\Vert^2_{V*} -\frac{\alpha}{2}\Vert\tilde\chi_t(t)\Vert^2
+ C_\alpha e^{-2t},
\end{equation}
for all $t\geq 0$.

Then, combining \eqref{unifest4} with \eqref{Lyap1}, we obtain
\begin{align}
\label{EQ1}
&\frac{d}{dt} \left( {\mathcal L} - \mu \langle {\mathbf q},\nabla A^{-1}_0\tilde\teta\rangle \right)
+ \Vert{\mathbf q}\Vert^2 + \Vert\tilde\chi_t\Vert^2_{V*} + \alpha\Vert\tilde\chi_t\Vert^2
+ \mu\Vert \tilde\teta\Vert^2
\\
\nonumber
&- \mu\langle {\mathbf q},\nabla A^{-1}_0\tilde\teta\rangle
- \mu\langle {\mathbf q},\nabla A^{-1}_0\tilde\chi_t\rangle
+ \mu \Vert A^{-1/2}_0\nabla\cdot{\mathbf q}\Vert^2 = \langle h(\tilde\chi),\tilde\chi_t\rangle ,
\end{align}
for some $\mu>0$ to be chosen below.

Following a well-known strategy (see, e.g., \cite{Je2,CJ}) we consider the functional
$$
{\mathcal G}(t) = \langle A_0^{-1}\tilde \chi_t ,A_0^{-1}(A_0\tilde\chi
+ \tilde\phi(\tilde \chi) - \overline{\tilde\phi(\tilde \chi)}) \rangle, \qquad t\geq 0,
$$
where
$$
\overline{\tilde\phi(\tilde \chi)}=<\tilde\phi(\chi),1>,
$$
and we observe that
\begin{align}
\label{EQ2}
\frac{d}{dt} {\mathcal G} &= \langle A_0^{-1}\tilde\chi_{tt},A_0^{-1}(A_0\tilde\chi + \tilde\phi(\tilde\chi) -
\overline{\tilde\phi(\tilde \chi)})\rangle\\
\nonumber
&+ \langle A_0^{-1}\tilde\chi_t,A_0^{-1}(A_0\tilde \chi_t + \tilde\phi'(\tilde \chi)\tilde \chi_t -
\overline{\tilde\phi'(\tilde \chi)\tilde \chi_t})\rangle\\
\nonumber
&= -\langle A_0^{-1}\tilde \chi_{t},A^{-1}_0(A_0\tilde\chi + \tilde\phi(\tilde\chi) -
\overline{\tilde\phi(\tilde \chi)}) \rangle\\
\nonumber
&-  \Vert A_0^{-1/2}(A_0\tilde\chi + \tilde\phi(\tilde\chi) - \overline{\tilde\phi(\tilde \chi)})\Vert^2\\
\nonumber
&-\alpha\langle \tilde \chi_{t},
A_0^{-1}(A_0\tilde\chi + \tilde\phi(\tilde\chi) - \overline{\tilde\phi(\tilde\chi)}) \rangle\\
\nonumber
&+ \langle \tilde \teta,
A_0^{-1}(A_0\tilde \chi + \tilde\phi(\tilde \chi) - \overline{\tilde\phi(\tilde\chi)}) \rangle\\
\nonumber
&+ \langle h(\tilde \chi),
A_0^{-1}(A_0\tilde \chi + \tilde\phi(\tilde \chi) - \overline{\tilde\phi(\tilde\chi)}) \rangle\\
\nonumber
&+ \Vert A^{-1/2}_0\tilde\chi_t\Vert^2
+ \langle A^{-1}_0\tilde\chi_t,A_0^{-1}(\tilde\phi'(\tilde\chi)\tilde\chi_t -
\overline{\tilde\phi'(\chi)\tilde \chi_t})\rangle.
\end{align}
Observe that (cf.~\eqref{fi2})
\begin{equation}
\label{EQ3}
\langle A_0^{-1}\tilde\chi_t,A_0^{-1}(\tilde\phi'(\tilde \chi)\tilde \chi_t -
\overline{\tilde\phi'(\tilde \chi)\tilde \chi_t})\rangle
\leq C\Vert A_0^{-1/2}\tilde \chi_t\Vert^2.
\end{equation}
Then, from \eqref{EQ1} and \eqref{EQ2}, using the Young inequality, we find (cf. also \eqref{newvar2},
\eqref{Lyap2}, and \eqref{EQ3})
\begin{equation}
\label{EQ4}
\frac{d}{dt}{\mathcal M} +  C_{\mu,\nu}{\mathcal N}^2\leq 0,
\end{equation}
for $\mu>0$ and $\nu>0$ sufficiently small, where
\begin{align}
\label{EQ5}
{\mathcal M}(t)&=\frac{1}{2}\left(\Vert\tilde\teta\Vert^2 + \Vert {\mathbf q}\Vert^2
+ \Vert\tilde\chi_t\Vert^2_{V^*}\right)+ E(\tilde\chi) - E(\tilde\chi_\infty)\\
\nonumber
&- \mu \langle {\mathbf q},\nabla A^{-1}_0\tilde\teta\rangle + \nu {\mathcal G} + C_{\alpha,\nu} e^{-2t},\\
\label{EQ6}
{\mathcal N}^2(t)&= \Vert{\mathbf q}(t)\Vert^2 + \Vert\tilde\chi_t(t)\Vert^2_{V^*} + \mu\Vert \tilde\teta(t)\Vert^2\\
\nonumber
&+ \nu\Vert A_0\tilde\chi(t) + \tilde\phi(\tilde\chi(t)) - \overline{\tilde\phi(\tilde \chi(t))} \Vert^2_{V^{-1}_0},
\end{align}
for all $t\geq 0$.

Let us introduce the unbounded set
$$
\Sigma = \big\{t\ge 0\,:\, \Vert \tilde \chi (t) - \tilde\chi_\infty\Vert_{V^1_0} \le \frac\eta 3\big\}
$$
where $\eta$ is given by Lemma~\ref{LSINE}. Then, for every $t\in\Sigma$, define
$$
\tau(t) = \sup \big\{t'\ge t\,:\, \sup_{s\in[t,t']} \Vert \tilde\chi(s) - \tilde\chi_\infty\Vert_{V^1_0}
\le \eta \big\},
$$
and observe that $\tau(t)>t$, for every $t\in\Sigma$.

Recalling \eqref{C1}-\eqref{C3}, let $t_0\in\Sigma$ be large enough such that
\begin{equation}
\label{EQ7}
\Vert\tilde\teta(t)\Vert +
\Vert{\mathbf q}(t)\Vert + \Vert \tilde \chi_t(t)\Vert_{V^*}
\le 1, \qquad \forall\,t\ge t_0,
\end{equation}
and set
\begin{align*}
J&=[t_0,\tau(t_0)),\\
J_1 &= \left\{t\in J\,:\, {\mathcal N}(t)
> e^{-2(1-\rho)t} \right\},\\
J_2 &=J\setminus J_1.
\end{align*}From \eqref{EQ4}, we have that ${\mathcal M}$ is decreasing, therefore it
is constant on $\omega(\teta_0,{\mathbf q}_0,\chi_0,\chi_1)$. In addition,
there holds
\begin{equation}
\label{EQ8}
\frac{d}{dt}\big(\vert {\mathcal M}(t)\vert^\rho\,{\rm sgn}\,{\mathcal
 M}(t)\big)
 = \rho\vert {\mathcal M}(t)\vert^{\rho-1}\frac{d}{dt} {\mathcal M}(t),\qquad \forall\,t\geq 0,
\end{equation}
so that $\vert {\mathcal M}\vert^\rho\,{\rm sgn}\,{\mathcal M}$ is decreasing as well.

Observe now that, for every $t\in J_1$, thanks to \eqref{LS} and \eqref{EQ7}, we have
$$
\vert {\mathcal M}(t)\vert^{1-\rho} \le
C{\mathcal N}(t),
$$
possibly choosing $\mu$ and $\nu$ even smaller than before.
Then, on account of \eqref{EQ4} and \eqref{EQ8}, we infer
\begin{align}
\label{EQ9}
\int_{J_1} {\mathcal N}(t) dt
&\le -C\int_{t_0}^{\tau(t_0)} \frac{d}{dt}\big(\vert {\mathcal M}(t)\vert^\rho\,
{\rm sgn}\,{\mathcal M}(t)\big)dt\\
\nonumber
&\le C\Big (\vert {\mathcal M}(t_0)\vert^\rho
+ \vert {\mathcal M}(\tau(t_0))\vert^\rho\Big),
\end{align}
where we mean that $\vert {\mathcal M}(\tau(t_0))\vert=0$ if $\tau(t_0)=\infty$.
On the other hand, we easily get
$$
\int_{J_2}{\mathcal N}(t)dt\leq Ce^{-2(1-\rho)t_0}.
$$
Therefore $\Vert \tilde\chi_t(\cdot)\Vert_{V^*}$ is integrable over $J$ and
\begin{align}
\label{EQ10}
0 &\le \limsup_{t_0\in\Sigma,\,t_0\to \infty} \,\int_{t_0}^{\tau(t_0)}
\Vert \tilde\chi_t(t)\Vert_{V^*} dt\\
\nonumber
&\le c\limsup_{t_0\in\Sigma,\,t_0\to \infty}\,\Big(\vert {\mathcal M}(t_0)\vert^\rho
+ \vert {\mathcal M}(\tau(t_0))\vert^\rho
+ Ce^{-2(1-\rho)t_0}\Big)=0.
\end{align}
Notice that, for every $t\in J$,
\begin{equation}
\label{EQ11}
\Vert \tilde\chi(t) - \tilde\chi_\infty \Vert_{V^*} \le \int_{t_0}^t \Vert \tilde\chi_t(s)\Vert_{V^*} ds
+ \Vert \tilde\chi (t_0) - \tilde\chi_\infty\Vert_{V^*}.
\end{equation}
Suppose now that $\tau(t_0)<\infty$ for any $t_0\in\Sigma$.
Then, by definition,
$$
\Vert \tilde\chi(\tau(t_0)) - \tilde\chi_\infty\Vert_{V^1_0} = \eta, \qquad \forall\,t_0\in\Sigma.
$$
Consider an unbounded sequence $\{t_n\}_{n\in\N}\subset\Sigma$ with the property
$$
\lim_{n\to \infty} \Vert \tilde\chi(t_n) - \tilde\chi_\infty\Vert_{V^1_0} =0.
$$
By compactness, we can find a subsequence $\{t_{n_k}\}_{k\in\N}$ and
an element $v_\infty\in \D(A)$ such that $(\teta_\infty,{\mathbf 0},v_\infty,0)\in
\omega(\teta_0,{\mathbf q}_0,\chi_0,\chi_1)$,
$\Vert \tilde v_\infty - \tilde\chi_\infty\Vert_{V^1_0} = \eta$, and
$$
\lim_{k\to \infty} \Vert \tilde\chi(\tau(t_{n_k})) -  \tilde v_\infty\Vert_{V^1_0} =0.
$$
Then, owing to \eqref{EQ10} and \eqref{EQ11}, we deduce the contradiction
$$
0<\Vert \tilde v_\infty - \tilde \chi_\infty \Vert_{V^*}\le \limsup_{k\to \infty}
\left(\int_{t_{n_k}}^{\tau(t_{n_k})} \Vert \tilde\chi_t(s)\Vert_{V^*} ds
+ \Vert \tilde\chi (t_{n_k}) - \tilde \chi_\infty\Vert_{V^*}\right) =0.
$$
Hence, $\tau(t_0)=\infty$ for some $t_0>0$ large enough and, recalling \eqref{MC1},
we can deduce that $\Vert \chi_t(\cdot)\Vert_{V^*}$ is indeed integrable over $(t_0,\infty)$.
This yields \eqref{convchi} by precompactness. On the other hand, on account of \eqref{C1}-\eqref{C3},
\eqref{singleton} holds as well. Finally, arguing as in \cite{GPS2}, we can prove that
\begin{equation}
\label{rate2}
\int_t^\infty {\mathcal N}(\tau)d\tau \leq C t^{-\frac{\rho}{1-2\rho}},\qquad\forall\,t\geq t^*,
\end{equation}
for some $t^*>0$.
This entails (cf.~\eqref{MC1} and \eqref{EQ6})
$$
\int_t^\infty \Vert \chi_t(\tau)\Vert_{V^*} d\tau \leq C t^{-\frac{\rho}{1-2\rho}},
\qquad\forall\,t\geq t^*.
$$
Thus we have
\begin{equation}
\label{rate3}
\|\chi(t) - \chi_\infty\|_{V^*}  \le C t^{-\frac{\rho}{1-2\rho}}, \qquad\forall\,t\geq t^*.
\end{equation}
Recalling now \eqref{MC1}, setting $v=1$ in the first equation
of \eqref{PHabs}, and integrating from $t$ to $\infty$, we obtain
\begin{equation}
\label{rate4}
\langle\teta(t),1\rangle - \teta_\infty = \langle \chi_1,1 \rangle e^{-t}.
\end{equation}
Therefore, integrating the first equation of
\eqref{PHabsav} with respect to time from $t\geq t^*$ to $\infty$, we deduce
$$
\langle -\tilde\teta(t) + \tilde\chi_\infty - \tilde \chi(t),v\rangle
= \int_t^\infty \langle {\mathbf q}(\tau),\nabla v\rangle d\tau,
$$
so that, on account of \eqref{EQ6}, \eqref{rate2}, and \eqref{rate3}, we infer
\begin{equation}
\label{rate5}
\|\tilde\teta(t)\|_{V^*}  \le C t^{-\frac{\rho}{1-2\rho}}, \qquad\forall\,t\geq t^*.
\end{equation}
Therefore, rate estimate \eqref{rate1} is a consequence of \eqref{rate3}-\eqref{rate5}.
In the case $\sigma=0$ we can proceed in a similar (actually, simpler) way,
noting that ${\mathbf q}=-\nabla\teta$.
\end{proof}
\begin{remark}
The decay estimate \eqref{rate1} for $\teta$
can be slightly improved. Actually,
using the decomposition 
\begin{equation}
\label{giulio1}
\Vert\teta(t) - \teta_\infty\Vert^2
\le 2\Vert\teta^d(t)\Vert^2+2\Vert\teta^c(t) -\teta_\infty\Vert^2,
\end{equation}
we see that, by \eqref{decay3}, the first term decays exponentially. 
Concerning the latter, one can use \eqref{smooth4}, \eqref{rate5}
and the interpolation inequality $\|v\|^2\le c\|v\|_{V}\|v\|_{V^*}$,
holding for all $v\in V$. Thus, \eqref{giulio1} eventually gives
\begin{equation}
\label{rate6}
  \|\teta(t)-\teta_\infty\|  \le C t^{-\frac{\rho}{2-4\rho}}, \qquad\forall\,t\geq t^*.
\end{equation}
\end{remark}


\section{Appendix}

\noindent

\newcommand{\bFormula}[1]{\begin{equation} \label{#1}}
\newcommand{\eF}{\end{equation}}

\newcommand{\vc}[1]{ {\bf #1} }
\font\F=msbm10
\newcommand{\dx}{dx}
\newcommand{\intO}[1]{ \int_{\Omega} #1 \,\dx }
\newcommand{\Grad}{\nabla}
\newcommand{\ov}{v_\infty}
\noindent
This section is devoted to demonstrate Lemma~\ref{LSINE}. Let us introduce the functional
\[
E(v) = \intO{ \left(\frac{1}{2} |\Grad v |^2 +  \tilde\Phi(v)\right) },
\]
defined on the space $V^1_0$. As before, we assume  $|\Omega|=1$.
The differential operator associated with the gradient $\partial
E$ does not conserve null mean functions. Hence the version of the
\L ojasiewicz theorem given in \cite{HJ} is not applicable
directly. This problem was solved in \cite{GG}, but
\cite[Assumption~5]{GG} is not exactly satisfied here. Our proof,
essentially given for the reader's convenience, follows the lines
of \cite{CFP} based on a general version of the \L
ojasiewicz-Simon theorem obtained in \cite{Chi}.

\smallskip

{\it Proof of Lemma 4.1.} We begin to observe that $\ov$
satisfying (\ref{S}) is a solution to
\bFormula{a1}
A\ov+\tilde\phi(\ov)- \overline{\tilde\phi(\ov)}=0.
\eF
Moreover,
$\ov$ is a critical point of $E$ on $V^1_0$. Indeed, it is easy to
check that, owing to our hypotheses, $E$ is continuously
differentiable on $V_0^1$, and
\begin{align}
\label{a2}
\partial E(\ov)h &= \intO{\left(
 \Grad \ov \cdot \Grad h + \tilde\phi(\ov)h \right)}\\
 &=\intO{\left(\Grad \ov \cdot \Grad h + \tilde\phi(\ov)h-\overline{ \tilde\phi(\ov)}h\right)},
 \nonumber
\end{align}
for all $h \in V_0^1$.

We recall that the dual space $(V_0^1)^*$ is the space of classes
\[
[f]=\{f+g;\ g\in V^*,<<g,V_0^1>>=0\}, \quad f\in V^*,
\]
where $<<\cdot,\cdot>>$ stands for the duality between $V^*$ and
$V$, endowed with the norm
\[
\|[f]\|_{(V_0^1)^*}=\inf_{g\in V^*,<g,V^1_0>=0}\|f+g\|_{V^*}
=\inf_{c\in \R}\|f+c\|_{V^*}=\|f-\overline f\|_{V^*} .
\]

Consider the mapping $F: V_0^2 \to H_0$, $F = \partial E|_{V_0^2}$ defined by
\[
F(v) =A_0v+\tilde\phi(v) -\overline{\tilde\phi(v)}.
\]
By virtue a well-known Sobolev embedding theorem, we have that
$v_\infty \in V^2_0\subset L^\infty(\Omega)$, and, due to our
assumptions, we can find a neighborhood ${\mathcal U}(\ov)$ in the
space $V_0^2$ such that $F: {\mathcal U}(\ov ) \to H_0$ is
analytic. Further, $A_0:\ V_0^1\to (V_0^1)^*$ and $A_0:\ V_0^2\to
H_0$ have compact resolvents. Observe now that, when $\partial^2 E
: V_0^1 \to {\mathcal Lin}[V_0^1,(V_0^1)^*]$, then
\bFormula{a3}
<<\partial^2 E (v)[w], z>>
= \intO{ \tilde\phi'(v) wz } + \intO{\Grad w \cdot \Grad z }.
\eF
In addition, we have
\[
\partial^2 E =\partial F :
 V_0^2 \to {\mathcal Lin}[V_0^2,H_0], \quad
  \partial F(v)[w]= A_0
 w+\tilde\phi(v)w-\overline{\tilde\phi(v)w}.
\]
Hence, in both cases, $\partial^2 E(v_\infty)$ can be viewed as a
bounded perturbation of $A_0$ restricted to the respective spaces.
It follows that $\mbox{Ker}\ \partial^2 E (v_\infty)\subset V_0^2$
and its range is closed in $(V_0^1)^*$ and $H_0$, respectively.
Moreover, there holds
\[
(V_0^1)^*=\mbox{Ker}\ (\partial^2 E(v_\infty))\oplus \mbox{Ran}\
(\partial^2 E(v_\infty)),\ \ H_0=\mbox{Ker}\ (\partial^2
E(v_\infty))\oplus \mbox{Ran} \ (\partial F(v_\infty)).
\]
Now, we can apply \cite[Thm.~3.10]{Chi} and \cite[Cor.~3.11]{Chi}
to obtain
$$
| E(v) - E(v_\infty)|^{1-\rho} \le L\|\partial E(v)
\|_{(V_0^1)^*},
$$
and, consequently, (\ref{LS}).


\end{document}